\documentstyle[12pt]{article}
\newcommand{\fr}{ \frac}
\newcommand{\beq}{\begin{equation}}
\newcommand{\eeq}{\end{equation}}
\newcommand{\br}{\begin{array}}
\newcommand{\er}{\end{array}}
\newcommand{\lb}{\label}
\newcommand{\r}{\ref}
\newcommand{\ar}{\rightarrow}

\newcommand{\ov}{\overline}
\newcommand{\lan}{\langle}
\newcommand{\rn}{\rangle}
\newcommand{\pr}{\prime}
\newcommand{\s}{\sigma}
\newcommand{\lm}{\lambda}

\begin{document}
\begin{center}
{\large \bf Unitary Representations of the 2-Dimensional Euclidean Group 
in the Heisenberg Algebra}
\end{center}

\vspace{1cm}

\begin{flushleft}
H. Ahmedov$^1$  and I. H. Duru$^{2,1}$
\vspace{.5cm}

{\small 
1. Feza G\"ursey Institute,  P.O. Box 6, 81220,  \c{C}engelk\"{o}y, 
Istanbul, Turkey 
\footnote{E--mail : hagi@gursey.gov.tr and duru@gursey.gov.tr}.

2. Trakya University, Mathematics Department, P.O. Box 126, 
Edirne, Turkey.}
\end{flushleft}
\vspace{1cm} \noindent 
{\bf Abstract}: $E(2)$ is studied  as the automorphism group 
of the Heisenberg algebra $H$. 
The basis in the Hilbert space $K$ of functions on $H$ on which 
the unitary irreducible representations of the group are realized 
is explicitly constructed. The addition theorem for the Kummer 
functions is derived. 

\begin{center}
Febriary  2000
\end{center}

\noindent
{\bf 1. Introduction}

\vspace{5mm}
\noindent
Investigating the properties of manifolds by means of the symmetries they 
admit has a long history. Non-commutative geometries have become the subject 
of similar studies in recent decades. For example there exists an extensive 
literature on the $q$-deformed groups $E_q(2)$ and $SU_q(2)$ which are the 
automorphism groups of the quantum plane $zz^* =qz^*z$ and the quantum sphere 
respectively \cite{w}. Using group theoretical methods the invariant distance 
and the Green functions have also been written in these deformed spaces
\cite{a}.

The purpose of the present work is to analyze yet another non-commutative 
space $[z, z^*]=\s$ ( i. e. the space generated by the Heisenberg  algebra )
by means of its automorphism group $E(2)$.

In Section 2 we define $E(2)$ in the Heisenberg algebra $H$ and construct 
the unitary representations of the group in the Hilbert space $X$ where 
$H$ is realized.

In Section 3 we give the unitary irreducible representations of $E(2)$ in 
the Hilbert space $K$ of the square integrable functions  on $H$ and 
construct the basis in $K$ where the irreducible representations of the group 
are realized. The basis are found to be written in terms of the Kummer 
functions. Commutative limit as $\s\ar 0$ is also discussed.

Section 4 is devoted to the addition theorem for the Kummer functions.
This theorem provides a group theoretical interpretation for the already 
existing identities involving the Kummer and Bessel functions. It may also 
lead to new identities. 
\newpage

\noindent
{\bf 2. $E(2)$ as the automorphism group  of the Heisenberg algebra}

\vspace{5mm}
\noindent
The one dimensional Heisenberg algebra $H$ is the 
3-dimensional vector space with the basis elements $\{ z, z^*, 1\}$ 
and the bilinear antisymmetric product  
\beq\lb{com}
[z, z^*]=1.
\eeq
The $*$-representation of $H$ in the suitable dense subspace of the Hilbert 
space $X$ with the complete orthonormal basis $\{\mid n\rn \}$, 
$n=0,1,2, ...$ is given by 
\beq\lb{r1}
z \mid n\rn =\sqrt{n} \mid n-1\rn, \ \ \ 
z^*\mid n\rn =\sqrt{n+1} \mid n+1\rn. 
\eeq
Let us represent the Euclidean group $E(2)$ in the vector space $H$ 
\beq\lb{tr}
g\left(
\begin{array}{c} 
z \\ 
z^* \\
1 
\end{array} 
\right) =
\left(
\begin{array}{ccc} 
e^{i \phi}  & 0 &  re^{i\psi} \\
0 & e^{-i \phi} &  re^{-i\psi} \\
0 & 0 &  1
\end{array} 
\right) 
\left(
\begin{array}{c} 
z \\ 
z^* \\
1 
\end{array} 
\right). 
\eeq
Since these transformations preserves the commutation relation
\beq
[g z, gz^*] =[z, z^*]
\eeq
we conclude that 
\beq\lb{u}
g z = U(g) z U^{-1}(g), \ \ \ g z^* = U(g) z^* U^{-1}(g)
\eeq
where $U(g)$ is the unitary representation of $E(2)$ in  $X$: 
\beq
U(g_1)U(g_2)=U(g_1g_2), \ \ \ U^*(g)=U^{-1}(g) = U(g^{-1}).
\eeq 
Simple verification shows that the normalized state
\beq
\mid 0 \rn^\pr = e^{-\fr{r^2}{2}}e^{re^{i(\psi-\phi)}z^* } \mid 0\rn
\eeq
satisfies the condition 
\beq
gz \mid 0\rn^\pr = 0.
\eeq
In the new orthonormal basis 
\beq
\mid n\rn^\pr = \fr{(gz^*)^n}{\sqrt{n!}} \mid 0\rn^\pr
\eeq
we have
\beq\lb{r2}
gz \mid n\rn^\pr =\sqrt{n} \mid n-1\rn^\pr, \ \ \ 
gz^* \mid n\rn^\pr =\sqrt{n+1} \mid n+1\rn^\pr. 
\eeq
The unitary operator $U(g)$ defines the transition between two  
complete orthonormal basis $\{\mid n\rn \}$ and $\{\mid n\rn^\pr \}$ 
\beq
\mid n\rn^\pr = U(g)\mid n\rn, \ \ \ 
\mid n\rn = U(g^{-1})\mid n\rn^\pr.
\eeq
Matrix elements of $U(g)$ in the basis $\{\mid n\rn \}$ reads
\beq
\lan m\mid U(g) \mid n\rn = 
\fr{e^{-\fr{r^2}{2}} }{\sqrt{n!}}
\lan m\mid (gz^*)^n e^{ - re^{i(\psi-\phi)}z^* } \mid 0\rn
\eeq
which after some simple algebra can be expressed in terms of the 
degenerate hypergeometric functions  as
\beq\lb{unit}
\lan m\mid U(g) \mid n\rn =  (-)^m e^{i(m-n)\psi - i n \phi} 
\fr{r^{n+m}e^{-\fr{r^2}{2}} }{\sqrt{n!m!}}\\_2F_0 (-m, -n; -\fr{1}{r^2}).
\eeq
Using the relations \cite{e1}
\beq
\\_2F_0 (-m, -n; -\fr{1}{r^2})=\fr{n!}{(n-m)!} (-\fr{1}{r^2})^m   
\Phi(-m,1+n-m; r^2), \ \ \ \ n\geq m,
\eeq
\beq
\\_2F_0 (-m, -n; -\fr{1}{r^2})=\fr{m!}{(m-n)!} (-\fr{1}{r^2})^n   
\Phi(-n,1+m-n; r^2), \ \ \ \ m\geq n
\eeq
we can also express the matrix elements in terms of the Kummer 
function $\Phi$.
 
\vspace{1cm}
\noindent

\noindent
{\bf 3. Unitary representations of $E(2)$ in the space of functions on $H$ }

\vspace{5mm}
\noindent
Let $K_0$ be set of finite sums
\beq\lb{set}
F = \sum (f_n (\zeta )z^n + z^{*n}f_{-n}(\zeta)). 
\eeq
Here  $f_n(\zeta)$ are functions of $\zeta =z^*z$ with finite support in 
$Spect (\zeta ) =\{0, 1, 2, ...\}$. Completion of  $K_0$ in the norm
\beq
\mid \mid F\mid\mid =\sqrt{tr (F^*F)}
\eeq
forms the Hilbert space  $K$ of the square integrable functions 
in the linear space $H$ with the scalar product 
\beq\lb{sc}
(F,G)=tr (F^*G).
\eeq
The formula 
\beq\lb{rep}
T(g) F(z) = F(gz)
\eeq
defines the representation of $E(2)$ in $K$. (\r{r1}) and (\r{r2}) and the 
independence of the trace from the basis over which it is taken  
imply that the representation is unitary. Using (\r{u}) we can rewrite 
(\r{rep}) in the form
\beq\lb{rep1}
T(g) F(z) = U(g)F(z) U^*(g). 
\eeq
Now we consider the infinitesimal form of (\r{rep}). 
Let $g=g(re^{i\psi}, \phi )$ in (\r{tr}). With the one parameter 
subgroups $g_1= g(\epsilon, 0)$, $g_2= g(i\epsilon, 0)$ and 
$g_3=g(0,\epsilon)$ of $E(2)$ we associate the linear operators
$K_0\ar K$
\beq\lb{inf}
p_k (F) = \lim_{\epsilon\ar 0} \fr{1}{\epsilon}(T(g_k)F - F) 
\eeq
with  the  limit taken in the strong operator topology.  
Inserting (\r{rep1}) into (\r{inf}) we get 
\beq\lb{real}
p (F) = 2[F, z^*] , \ \ \ \ov{p} (F) = 2[z,F], \ \ \
h (F) = [\zeta, F] ,  
\eeq
where 
\beq
p=p_1-ip_2, \ \ \ \ov{p}= p_1+ip_2, \ \ \ h=i p_3. 
\eeq
For example
\beq
p (f(\zeta)z^n)= 2( nf(\zeta) +\zeta ( f(\zeta +1) -f(\zeta)) ) z^{n-1}.
\eeq
(\r{sc}) and (\r{real}) imply the real structure in the Lie algebra 
of $E(2)$
\beq
p^*=-\ov{p}, \ \ \ h^*=h.
\eeq
The irreducible representations of $E(2)$ defined by the weight  
$\lm\in R$ can be constructed in the space of  square integrable 
functions on the circle and the matrix elements are given in terms of 
the Bessel functions \cite{v} 
\beq\lb{mat}
t^\lm_{kn}(g) = i^{n-k}e^{-i(n\phi+(k-n)\psi)}J_{n-k}(\lm r).
\eeq
Coming to our case the basis $D_k^\lm$ in $K$ where the unitary irreducible 
representations of the group are realized will be the eigenfunctions of 
the complete set of commuting operators $pp^*$ and $h$:
\beq\lb{c1}
pp^* D_k^\lm =\lm^2 D_k^\lm, 
\eeq
\beq\lb{c2}
h D_k^\lm =  k D_k^\lm .
\eeq
The solutions of the equation (\r{c2}) are 
\beq
D_k^\lm (z)= 
\{
\begin{array}{cc}
z^{*k}f_k^\lm (\zeta) &   \ \ if \ \ k\geq 0 \\ 
f_{-k}^\lm (\zeta) z^{-k}  &   \ \ if \ \ k\leq 0
\end{array}
\eeq
Inserting (\r{c2}) into (\r{c1}) we get
\beq
(k+1+\zeta )f^\lm_k(\zeta+1) + 
(\fr{\lm^2}{4}-2\zeta -k-1)f^\lm_k(\zeta)+\zeta f^\lm_k(\zeta -1) = 0. 
\eeq
By the virtue of the recurrence relation  \cite{e1}
\beq
a\Phi (a+1, b; c) + (a-b)\Phi (a-1, b; c) +(b-2a-c)\Phi (a, b; c)=0
\eeq
we observe that  the solutions are given in terms of the Kummer functions 
as
\beq
f^\lm_k (\zeta ) = \fr{(-\lm^2)^{k/2}}{2^k k!} e^{-\fr{\lm^2}{8}}
\Phi (-\zeta, 1+k; \fr{\lm^2}{4}). 
\eeq
The formula 
\beq
L^k_n (x) =\fr{(k+n)!}{k!n!} \Phi (-n, 1+k; x)
\eeq
allows us to express the basis elements  in terms of the Laguerre
polynomials too
\beq
f^\lm_k (\zeta ) = \fr{(-\lm^2)^{k/2} \zeta ! }{2^k (k+\zeta)!} 
e^{-\fr{\lm^2}{8}} L^k_\zeta (\fr{\lm^2}{4}).
\eeq
The above formula is well defined since the spectrum of the operator $\zeta$ 
is the set of positive integers with zero. It has been well known that the 
Laguerre polynomials are related to the group of complex third order 
triangular matrices  \cite{v}. The group parameters appears in the 
argument of the Laguerre polynomials. 
However in our case the group parameter appear in  the index of 
this function. 

To obtain the orthogonality relations we first take $z\ar 1^-$ limit in the 
summation formula \cite{e}
\beq
\sum_{n=0}^\infty \fr{n!}{(n+k)!}L^k_n(x)L^k_n(y)z^n= 
\fr{(xyz)^{-k/2}}{1-z}e^{-z\fr{x+y}{1-z}} I_k (2\fr{\sqrt{xyz}}{1-z}), \ \ 
\mid z\mid <1,
\eeq
then use the asymptotic relation  for the Bessel functions \cite{r}
\beq
I_\nu (x) \sim \fr{e^x}{\sqrt{2\pi x}} + 
\fr{e^{-x+(\nu+1/2)\pi}}{\sqrt{2\pi x}}
\eeq
and employ the representation
\beq
\delta (x)= 
\lim_{\epsilon\ar 0}\fr{1}{\sqrt{\pi\epsilon}} e^{-\fr{x^2}{\epsilon}} 
\eeq
for the Dirac delta function. As the result we have 
\beq
(D^\lm_k, D^{\lm^\pr}_n)=\delta_{kn}\delta (\lm^2-\lm^{\pr 2}).
\eeq

Let us make the substitution $z\ar \fr{1}{\sqrt{\s}} z$ in (\r{com}). 
In $\s\ar 0$ limit the linear space $H$  reduces to the complex plane;
and $D^\lm_k(z)$  reduces to $t^\lm_{k0}(g)$ which are the restriction 
of the matrix elements (\r{mat}) on the complex plane $E(2)/U(1)$:
\beq
\lim_{\s\ar 0} D^{\sqrt{\s}\lm}_k(\fr{1}{\sqrt{\s}} re^{i\psi}) = t^\lm_{k0}(g).
\eeq
To prove the above limit we used the asymptotic  relation \cite{e1}
\beq
\lim_{a\ar \infty} \Phi (a,b; \fr{c}{a}) = 
\Gamma (b) \sqrt{c^{1-b}} I_{b-1} (2\sqrt{c}). 
\eeq

\vspace{1cm}
\noindent

\noindent
{\bf 4. The addition theorem}

\vspace{5mm}
\noindent
Let us consider the  representations of $E(2)$ in the 
basis $D^\lm_{k}$: 
\beq
T(g) D^\lm_{k}(z) = \sum_{n=-\infty}^\infty t_{kn}^\lm (g) D^\lm_n(z).
\eeq
By making use of (\r{rep1}) the above formula can be rewritten as
\beq\lb{add}
U(g) D^\lm_{k} (z) U^*(g) = 
\sum_{n=-\infty}^\infty t_{kn}^\lm (g) D^\lm_n (z)
\eeq
which is an addition theorem  useful in deriving  many  identities 
involving the Kummer and Bessel functions.
For example let $g=g(r,o)$, $k\geq  0$ and $x=\lm /2$. 
Then (\r{add}) reads
\begin{eqnarray}\lb{add1}
\fr{x^k}{k!} U z^{*k}\Phi (-\zeta , 1+k; x^2) U^* = 
\sum_{n=0}^\infty
\fr{x^n}{n!} J_{k-n} (2rx) z^{*n}\Phi (-\zeta , 1+n; x^2) + \nonumber \\
\sum_{n=1}^\infty
\fr{(-x)^n}{n!} J_{k+n} (2rx) \Phi (-\zeta , 1+n; x^2)z^n,
\end{eqnarray}
where
\beq
U = e^{-\fr{r^2}{2}} \sum_{n,m=0}^\infty  
\fr{(-)^mr^{n+m} }{\sqrt{n!m!}} \\_2F_0 (-m, -n; -\fr{1}{r^2})
\mid m\rn\lan n\mid.
\eeq
Sandwiching (\r{add1}) between the states $\lan 0\mid$ and $\mid 0\rn$ 
we get
\beq
\sum_{n=0}^\infty \fr{r^{2n}}{n!}\Phi (-n, 1+k; x^2)= k!
(xr)^{-k}e^{r^2} J_k (2xr).
\eeq
Multiplying (\r{add1}) by $U^*$ from the left and then sandwiching it 
between the states $\lan m+k\mid$ and $\mid 0\rn$  we obtain another 
formula 
\beq
\fr{(m+k)!}{m!k!} (\fr{x}{r})^k \Phi (-m, 1+k; x^2)=
\sum_{n=0}^\infty \fr{(-xr)^n}{n!} \\_2F_0(-m-k,-n; -\fr{1}{r^2})
J_{k-n} (2xr).
\eeq
It is clear that (\r{add}) can lead to many more identities, that some of
them may not exist in the literature.

\end{document}